\documentclass[a4paper,12pt]{article}
\usepackage{amsmath,amssymb,amsthm}
\usepackage{amscd}

\usepackage[all]{xy}

\makeatletter
\@addtoreset{equation}{subsection}
\makeatother

\begin{document}

\title{
A regular homotopy version of the Goldman-Turaev Lie bialgebra, 
the Enomoto-Satoh traces and the divergence cocycle 
in the Kashiwara-Vergne problem
}
\author{Nariya Kawazumi 
}
\date{May 31, 2014}
\maketitle

\abstract{
By introducing a refinement of the Goldman-Turaev Lie bialgebra, we interpret the divergence cocycle in the Kashiwara-Vergne problem and the Enomoto-Satoh obstructions for the surjectivity of the Johnson homomorphisms as some part of a regular homotopy version of the Turaev cobracket. 
}

\begin{center}
{\bf Introduction}
\end{center}
This is an announcement on my research in progress, 
which introduces a refinement of 
the Goldman-Turaev Lie bialgebra. 
Its goal is to interpret 
the divergence cocycle in the Kashiwara-Vergne problem \cite{AT} and 
the Enomoto-Satoh obstructions for the surjectivity 
of the Johnson homomorphisms (= the Enomoto-Satoh traces) \cite{ES}
as some part of a regular homotopy version of the Turaev cobracket. \par
In my previous work joint with Yusuke Kuno \cite{KK3} 
we proved that the Morita traces are included in the Turaev cobracket. 
The Enomoto-Satoh traces \cite{ES} are refinements of the Morita traces, 
and closely related to 
the divergence cocycle in the Kashiwara-Vergne problem. 
Enomoto \cite{E} proved that the graded quotient of the Turaev cobracket 
does {\it not} include the Enomoto-Satoh traces. 
This fact seems to come from the fact that the Turaev cobracket is defined to 
be invariant under the birth-death move of a monogon in free loops. 
This is the reason why we consider the {\it regular} homotopy set 
of {\it immersed} free loops on a surface. 
On the other hand, the first term of the Enomoto-Satoh traces 
is just the Earle class $k$ on the mapping class group.
Furuta gave an explicit cocycle for the Earle class in terms 
of a framing of the tangent bundle of the surface. 
For details, see \cite{Mo} \S4. Our construction is inspired by
Furuta's construction. \par
Proofs and details of these results will appear elsewhere.
The author thanks Naoya Enomoto, Takao Satoh and Yusuke Kuno 
for valuable discussions, and he
is partially supported by the Grant-in-Aid for
Scientific Research (S) (No.24224002) and (B) (No.24340010) from the
Japan Society for Promotion of Sciences.

\begin{center}
{\bf A regular homotopy version of the Goldman-Turaev Lie bialgebra}
\end{center}
Let $S$ be a compact connected oriented $C^\infty$ surface 
with $\partial S \neq \emptyset$. We denote by 
$\hat \pi^+= \hat\pi^+(S)$ the {\it regular} homotopy 
set of free {\it immersed} loops on $S$. The infinite cyclic group 
$\langle r \rangle$ acts on the set $\widehat{\pi}^+$ by inserting 
a (positive) monogon into a loop. The action is free, and the 
orbit space $\hat \pi = \hat\pi(S) := \hat\pi^+(S)/\langle r \rangle$ 
equals the free homotopy set of free loops on $S$. 
We denote by $\Phi: \hat\pi^+ \to \hat\pi$ the quotient map, which can be 
regarded as the map forgetting smooth structures on immersed loops. 
The rational group ring $\mathbb{Q}\langle r\rangle$ is naturally identified with
the Laurent polynomial ring $\mathbb{Q}[r, r^{-1}]$. 
The $\mathbb{Q}$-free vector space over the set $\hat \pi^+$,
$\mathbb{Q}\hat\pi^+$, is a free $\mathbb{Q}\langle r\rangle$-module. 
We denote by $\mathbb{Q}\langle r\rangle1$ 
the linear span of the regular homotopy classes of null-homotopic immersed 
free loops on $S$. 
\par
Since $S$ is connected and its boundary is non-empty,  
the tangent bundle $TS$ is trivial. 
We call the homotopy class $f$ of a global trivialization 
$TS \cong S\times\mathbb{R}^2 \overset{\operatorname{pr}_2}\to 
\mathbb{R}^2$ a {\it framing} of $S$. 
If we fix a framing $f$, then we can define the (global) rotation number
$\operatorname{rot}_f: \hat\pi^+ \to \mathbb{Z}$. 
The map $\tilde \Phi_f:= (\Phi, \operatorname{rot}_f): \hat\pi^+ \to \hat\pi\times
\mathbb{Z}$ is a bijection. We define $s_f: \hat\pi \to \hat\pi^+$ by 
$s_f(\alpha) := \tilde {\Phi_f}^{-1}(\alpha, 0)$ for $\alpha \in \hat\pi$, and 
$\varepsilon_f: \mathbb{Q}\hat\pi^+ \to \mathbb{Q}\langle r\rangle$ by 
$\varepsilon_f(\beta) := r^{\operatorname{rot}_f(\beta)}$ for $\beta \in \hat\pi^+$.
\par
The {\it regular} Goldman bracket $[\,, \,]^+: \mathbb{Q}\hat\pi^+ 
\otimes_{\mathbb{Q}\langle r\rangle}\mathbb{Q}\hat\pi^+ \to 
\mathbb{Q}\hat\pi^+$ is defined in the same way as the original one \cite{Go}. 
The {\it regular} Turaev cobracket 
$
\delta^+: \mathbb{Q}\hat\pi^+
\to \mathbb{Q}\hat\pi^+
\otimes_{\mathbb{Q}\langle r\rangle}
\mathbb{Q} \hat\pi^+
$
is also defined in a similar way to the original one \cite{T}. 
The triple $(\mathbb{Q}\hat\pi^+, [\,, \,]^+, \delta^+)$ is a Lie bialgebra.  
For any embedded loop $\alpha \in \hat\pi^+$ and $n \in \mathbb{Z}$ 
we have $\delta^+(\alpha^n) 
= 0$. In particular, the cobracket $\delta^+$ 
vanishes on $\mathbb{Q}\langle r\rangle 1$. Hence we obtain the induced 
operation
$$
\delta^+: \mathbb{Q}\hat\pi^+
/\mathbb{Q}\langle r\rangle 1
\to \mathbb{Q}\hat\pi^+
\otimes_{\mathbb{Q}\langle r\rangle}
\mathbb{Q} \hat\pi^+.
$$
In the original case \cite{T} the target of $\delta$ is 
$(\mathbb{Q}\hat\pi/\mathbb{Q}1)^{\otimes 2}$, since the cobracket has to 
be invariant under the birth-death move of a monogon, which we can ignore 
in the context of regular homotopy. \par
We number the connected components of the boundary $\partial S
= \coprod_{a=0}^n\partial_aS$, where $n = \sharp\pi_0(\partial S) - 1$. 
For each $a$ we choose a point $*_a \in 
\partial_aS$ and an inward vector $v_a \in T_{*_a}S$. 
We define a $\mathbb{Q}$-linear small category $\mathbb{Q}\Pi^+S\vert_E$ 
whose object set is $E := \{*_a\}^n_{a=0}$, and whose morphism vector space 
from $*_a$ to $*_b$ is the $\mathbb{Q}$-free vector space over the set 
$\Pi^+S(v_a, -v_b) := \{\ell: [0, 1] \to S; \, \text{an immersed path in $S$ 
from $*_a$ to $*_b$ with $\overset\centerdot\ell(0) = v_a$}\linebreak
\text{and 
$\overset\centerdot\ell(1) = -v_b$}\}$ 
modulo regular homotopy. 
The infinite cyclic group $\langle r\rangle$ acts on the set $\Pi^+S(v_a, -v_b)$
by inserting monogons into paths. 
If we fix a framing $f$ of $S$, we have a group isomorphism
$\Pi^+S(v_a, -v_a) \cong \pi_1(S, *_a)\times \mathbb{Z}$. 
We denote by $\operatorname{Der}_\partial
(\mathbb{Q}\Pi^+ S\vert_E)$ the Lie algebra of 
$\mathbb{Q}\langle r\rangle$-linear derivations of the category 
$\mathbb{Q}\Pi^+S\vert_E$ annihilating all loops parallel 
to some boundary component. 
In the same way as in \cite{KK1} 
we can define a $\mathbb{Q}\langle r\rangle$-Lie algebra homomorphism
$\sigma^+: \mathbb{Q}\hat\pi^+/\mathbb{Q}\langle r\rangle1\to 
\operatorname{Der}_\partial(\mathbb{Q}\Pi^+ S\vert_E)$. \par
Now we take completions of $\mathbb{Q}\hat\pi^+$ and 
$\mathbb{Q}\Pi^+ S\vert_E$ with respect to the augmentation ideal 
of the group ring $\mathbb{Q}\Pi^+S(v_a, -v_a)$, and denote them by 
$\widehat{\mathbb{Q}\hat\pi^+}$ and 
$\widehat{\mathbb{Q}\Pi^+ S}\vert_E$, respectively. 
Recall the completed group ring $\widehat{\mathbb{Q}\langle r\rangle}$ 
is naturally identified with the ring of formal power series in $\rho := \log r$. 
In other words, we have $\widehat{\mathbb{Q}\langle r\rangle} = 
\mathbb{Q}[[\log r]] = \mathbb{Q}[[\rho]]$. The bracket $[\,, \,]^+$ and 
the cobracket $\delta^+$ induce a natural Lie bialgebra structure on 
the completion $\widehat{\mathbb{Q}\hat\pi^+}$. \par

\begin{center}
{\bf The Enomoto-Satoh traces}
\end{center}
Recall that the completed Goldman-Turaev Lie bialgebra
$\widehat{\mathbb{Q}\hat\pi}$ introduced in \cite{KK2} 
has a decreasing filtration 
$\{\widehat{\mathbb{Q}\hat\pi}(m)\}_{m=1}^\infty$, and 
that the $\mathbb{Q}$-linear category $\widehat{\mathbb{Q}\Pi S}\vert_E$ 
admits a coproduct $\Delta$ \cite{KK2}.
Then we introduce a Lie subalgebra 
$
L^+(S) := \{u \in \widehat{\mathbb{Q}\hat\pi}(3); \,
(\sigma(u)\widehat{\otimes}1+1\widehat{\otimes}\sigma(u))\Delta = 
\Delta\sigma(u)\}\subset \widehat{\mathbb{Q}\hat\pi}.
$
We can prove that the restriction of the map 
$s_f: \widehat{\mathbb{Q}\hat\pi}
\to \widehat{\mathbb{Q}\hat\pi^+}/\mathbb{Q}[[\rho]]1$ 
to the subalgebra $L^+(S)$ 
does not depend on the choice of a framing $f$. So we denote it by 
$s_{\text{can}}: L^+(S) \to 
\widehat{\mathbb{Q}\hat\pi^+}/\mathbb{Q}[[\rho]]1$, and call it the 
{\it canonical section}.  Then we define the maps 
${\sf ES}^+_f:  \widehat{\mathbb{Q} \hat\pi^+}/\mathbb{Q}[[\rho]]1 \to  \widehat{\mathbb{Q} \hat\pi^+}/\mathbb{Q}[[\rho]]1$
and 
${\sf ES}_f: L^+(S) \to \widehat{\mathbb{Q}\hat\pi}$ 
by the following commutative diagram
$$
\begin{xymatrix}{
\widehat{\mathbb{Q} \hat\pi^+}/\mathbb{Q}[[\rho]]1
\ar[rr]^{\delta^+}\ar[rrd]_{{\sf ES}^+_f}&&
\widehat{\mathbb{Q} \hat\pi^+}\widehat{\otimes}_{\mathbb{Q}[[\rho]]}
\widehat{\mathbb{Q} \hat\pi^+} \ar[d]^{\varepsilon_f\widehat{\otimes}1}\\
L^+(S) \ar[u]^{s_{\text{can}}}\ar[rrd]_{{\sf ES}_f}
&&  \widehat{\mathbb{Q}\widehat{\pi}^+}/\mathbb{Q}[[\rho]]1 \ar[d]^{\Phi}\\
&& \widehat{\mathbb{Q}\hat\pi}.
}\end{xymatrix}
$$
In the case the boundary $\partial S$ is connected,
the graded quotient of the map ${\sf ES}_f$
$$
\operatorname{gr}({\sf ES}_f): \operatorname{gr}(L^+(S)) \to 
\operatorname{gr}(\widehat{\mathbb{Q}\hat\pi})
$$ 
is exactly the Enomoto-Satoh traces. 
On the other hand, if $S$ is of genus $0$, the Lie algebra 
$L^+(S)$ is isomorphic to an extension of the positive
part of the special derivation algebra $\mathfrak{sder}_n$, 
and $\widehat{\mathbb{Q}\hat\pi}$ to the space $\mathfrak{tr}_n$ 
in \cite{AT}. 
Then the graded quotient $\operatorname{gr}({\sf ES}_f)$ 
equals the restriction of the divergence cocycle $\operatorname{div}$ 
in the Kashiwara-Vergne problem \cite{AT}.
The proof of these facts is based on a tensorial description of the homotopy 
intersection form by Massuyeau and Turaev \cite{MT}. 
Hence the Enomoto-Satoh traces and the divergence cocycle 
are interpreted as 
some part of the regular Turaev cobracket. 

\begin{center}
{\bf The mapping class group}
\end{center}
The homomorphism 
$\sigma^+: \mathbb{Q}\hat\pi^+/\mathbb{Q}\langle r\rangle1\to 
\operatorname{Der}_\partial(\mathbb{Q}\Pi^+ S\vert_E)$
induces a $\mathbb{Q}[[\rho]]$-Lie algebra homomorphism
$
\sigma^+: \widehat{\mathbb{Q}\hat\pi^+}/\mathbb{Q}[[\rho]]1 \to \operatorname{Der}_\partial
(\widehat{\mathbb{Q}\Pi^+ S}\vert_E).
$
Then it is a Lie algebra isomomorphism 
\begin{equation}
\sigma^+: \widehat{\mathbb{Q}\hat\pi^+}/\mathbb{Q}[[\rho]]1 \overset\cong\to \operatorname{Der}_\partial
(\widehat{\mathbb{Q}\Pi^+ S}\vert_E).
\label{a1}\end{equation}
Moreover, for any framing $f$ of $S$, the map $\tilde\Phi_f$ 
induces an isomorphism 
\begin{equation}
\tilde\Phi_f: 
\widehat{\mathbb{Q}\hat\pi^+}/
\mathbb{Q}[[\rho]]1
\overset\cong\to  
\widehat{\mathbb{Q}\hat\pi}\widehat{\otimes} \mathbb{Q}[[\rho]].
\label{a2}\end{equation}
Let $\mathcal{I}^L(S)$ be the largest Torelli group in the sense of Putman 
\cite{P}. By the isomorphism (\ref{a1}) we can define the {\it geometric} Johnson 
homomorphism
$$
\tau^+: \mathcal{I}^L(S) \to \widehat{\mathbb{Q}\hat\pi^+}/
\mathbb{Q}[[\rho]]1
$$
in the same way as in \cite{KK2}. 
Applying 
a regular homotopy version of the logarithm formula for Dehn twists
\cite{KK1}\cite{KK2}\cite{MT}
to Putman's generators of $\mathcal{I}^L(S)$ \cite{P}, we can prove 
\begin{equation}
\delta^+\circ\tau^+ = 0: \mathcal{I}^L(S) \to 
 \widehat{\mathbb{Q} \hat\pi^+}\widehat{\otimes}_{\mathbb{Q}[[\rho]]}
\widehat{\mathbb{Q} \hat\pi^+}.
\end{equation}
In fact, the cobracket $\delta^+$ vanishes at 
any power of any embedded loop. Let $C$ be a simple closed curve in 
the interior of $S$ with $\pm [C] = 0 \in H_1(S; \mathbb{Z})$.
In other words, $C$ is a bounding simple closed curve. Then the Dehn 
twist $t_C$ along $C$ satisfies $\tau(t_C) \in L^+(S)$. 
Choose a framing $f$ of $S$. Let $h$ be the genus of the subsurface 
bounded by $C$. Then, under the isomorphism (\ref{a2}), 
the regular homotopy version of the 
logarithm formula says
$$
\tilde\Phi_f(\tau^+(t_C)) = \frac12 \vert (\log(C) + (1-2h)\rho)^2\vert
= \frac12\vert (1-2h)\rho\log(C) + (\log(C))^2\vert
= \frac12\vert (\log(C))^2\vert,
$$
since $C$ is null-homologous. 
We define the subgroup $\mathcal{K}(S) \subset \mathcal{I}^L(S)$ by 
those generated by such Dehn twists. In view of a theorem of Johnson
\cite{J}, $\mathcal{K}(S)$ is the Johnson kernel if the boundary of $S$
is connected. \par
Consequently, for any $S$, we obtain the commutative diagram
$$
\begin{xymatrix}{
\mathcal{I}^L(S) \ar[r]^{\tau^+}& \widehat{\mathbb{Q} \hat\pi^+}
/\mathbb{Q}[[\rho]]1
\ar[rr]^{\delta^+}
\ar[rrd]_{{\sf ES}^+_f}&&
\widehat{\mathbb{Q} \hat\pi^+}\widehat{\otimes}_{\mathbb{Q}[[\rho]]}
\widehat{\mathbb{Q} \hat\pi^+} \ar[d]^{\varepsilon_f\widehat{\otimes}1}\\
\mathcal{K}(S) \ar[u]^{\text{incl.}} \ar[r]^{\tau}& L^+(S) \ar[u]^{s_{\text{can}}}
\ar[rrd]_{{\sf ES}_f}
&&  \widehat{\mathbb{Q}\widehat{\pi}^+}/\mathbb{Q}[[\rho]]1 \ar[d]^{\Phi}\\
& && \widehat{\mathbb{Q}\hat\pi}.
}\end{xymatrix}
$$
In particular, we have ${\sf ES}_f\circ\tau\vert_{\mathcal{K}(S)} = 0: 
\mathcal{K}(S)
\to \widehat{\mathbb{Q}\hat\pi}$. This gives a geometric proof for the fact that the 
Enomoto-Satoh traces are obstructions for the surjectivity of the Johnson 
homomorphims. \par

\begin{center}
{\bf The genus $0$ case}
\end{center}
Let $n \geq 2$ be an integer. Here we study the framed pure braid group 
$FP_n$ on 
$n$ strands on the $2$-disk. This is nothing but the largest Torelli group of 
$S := \Sigma_{0, n+1}$. 
By capping each of the boundary components except one
by the surface $\Sigma_{1,1}$ we obtain an embedding of the surface 
$\iota: S\hookrightarrow \hat S:= \Sigma_{n, 1}$, which induces 
a Lie algebra homomorphism $\iota: \widehat{\mathbb{Q}\hat\pi}(S) 
\to \widehat{\mathbb{Q}\hat\pi}(\hat S)$.  
The pull-back of the Enomoto-Satoh trace $\operatorname{gr}({\sf ES}_f)$ 
by the map $\iota$ equals 
the divergence cocycle $\operatorname{div}$ in the Kashiwara-Vergne 
problem \cite{AT}
up to some low degree map $H_1(S)^{\otimes 2} \to H_1(S)$. 
Here $f$ is any framing of $\hat S$. 
In this case we have $\mathcal{K}(S) = \{1\}$. Instead we consider 
the commutator subgroup $[FP_n, FP_n]$. Choose a framing $f$ of $S$. 
For any simple closed curve $C_i$, $i = 1,2$, in $S$, we have 
$$
\tilde\Phi_f(\tau^+(t_{C_i})) 
= \frac{1}{2}\operatorname{rot}_f(C_i)\rho\vert \log(C_i)\vert 
+ \frac{1}{2}\vert (\log(C_i))^2\vert.
$$
Since the genus of $S$ is zero, the homology group $H_1(S; \mathbb{Q})$ 
is spanned by the boundary loops. In particular, $\vert\log(C_i)\vert$ is in the 
center of $\widehat{\mathbb{Q}\hat\pi}$. Hence we have 
$$
\left[\tilde\Phi_f(\tau^+(t_{C_1})), \tilde\Phi_f(\tau^+(t_{C_2}))\right]
= \left[\frac{1}{2}\vert (\log(C_1))^2\vert, \frac{1}{2}\vert (\log(C_2))^2\vert\right], 
$$
which is independent of the choice of the framing $f$. 
Thus we have $\tau^+\vert_{[FP_n, FP_n]} = s_{\text{can}}\circ \tau: 
[FP_n, FP_n] \to \widehat{\mathbb{Q} \hat\pi^+}/\mathbb{Q}[[\rho]]1$, and 
$({\sf ES}_f\circ\tau)(\varphi) = 0$ for any $\varphi \in [FP_n, FP_n]$.
\par

\vskip 10mm
\noindent
Department of Mathematical Sciences, \\
University of Tokyo \\
3-8-1 Komaba, Meguro-ku, Tokyo, \\
153-8914, JAPAN. \\
kawazumi@ms.u-tokyo.ac.jp\\
\\
\end{document}